\documentclass[12pt,reqno]{amsart}
\usepackage{geometry}                
\usepackage{graphicx}
\usepackage{amssymb}
\usepackage{epstopdf}
\usepackage{verbatim}
\usepackage{fullpage}
\usepackage{mathrsfs}
\def\Z{\mathbb Z}
\def\Q{\mathbb Q}
\def\F{\mathbb F}
\def\gf{\mathbb F}
\newcommand{\Aut}{\mathrm{Aut}}
\newcommand{\prob}{\mathrm{Prob}}
\newcommand{\GL}{\mathrm{GL}}
\newcommand{\isom}{\cong}
\newcommand{\dnd}{\nmid}
\newcommand{\C}{\mathscr C}
\newcommand{\B}{\mathscr{C}}

\newcommand{\leg}[2]{\left(\frac{#1}{#2}\right)}
\newcommand{\floor}[1]{\left\lfloor#1\right\rfloor}

\newtheorem{thm}{Theorem}
\newtheorem{cor}[thm]{Corollary}
\newtheorem{prop}[thm]{Proposition}
\newtheorem{lma}[thm]{Lemma}
\newtheorem{conj}[thm]{Conjecture}
\theoremstyle{remark}
\newtheorem{rmk}[thm]{Remark}

\title{A Cohen-Lenstra phenomenon for elliptic curves}
\author{Chantal David}
\address{Department of Mathematics and Statistics,
Concordia University,
1455 de Maisonneuve West,
Montreal, QC H3G 1M8, Canada} \email{cdavid@mathstat.concordia.ca}
\author{Ethan Smith}
\address{Department of Mathematics\\Liberty University\\1971 University Blvd\\MSC Box 710052\\ Lynchburg, VA 24502\\USA}
\email{ecsmith13@liberty.edu}

\begin{document}



\maketitle

\begin{abstract}
Given an elliptic curve $E$ and a finite Abelian group $G$, we consider the problem of counting
the number of primes $p$ for which the group of points modulo $p$ is isomorphic to $G$.
Under a certain conjecture concerning the distribution of primes in short intervals, we obtain an
asymptotic formula for this problem on average over a family of elliptic curves.
\end{abstract}

\section{Introduction.}

Let $E$ be an elliptic curve defined over the rational field $\Q$.
Given a prime $p$ where $E$ has good reduction, we consider the reduced curve, which we denote by $E_p$.
In previous work~\cite{DS:2013}, we studied the arithmetic function
\begin{equation*}
M_E(N):=\#\{p: \#E_p(\gf_p)=N\},
\end{equation*}
where $\#E_p(\gf_p)$ denotes the number of $\F_p$-rational points on the reduction.
Despite being such a natural object to study, it seems that the recent paper of Kowalski~\cite{Kow:2006} is the
first to introduce the function $M_E(N)$ and ask interesting questions about its behavior.
Kowalski's motivation for studying $M_E(N)$ stems from its relation with what he calls ``elliptic twins."
A prime $p$ is called an $E$-twin prime if there exists a prime $q\ne p$ such that $\#E_p(\F_p) = \#E_p(\F_q)$.
Then $p$ is an $E$-twin prime if and only if $M_E(\# E_p(\F_p)) > 1$.

The Hasse bound states that $\#E_p(\F_p)$ is never very far from $p+1$.  In particular,
$\left|p+1-\#E_p(\F_p)\right|<2\sqrt p$.
It follows that if $\#E_p(\F_p)=N$, then
\begin{equation}\label{p range}
N^-:=(\sqrt N-1)^2<p<(\sqrt N+1)^2=:N^+.
\end{equation}
Hence $M_E(N)$ is a finite number, satisfying the trivial bound
$M_E(N)\ll\sqrt N/\log(N+1)$. In \cite{Kow:2006}, the author shows that if $E$ possesses complex multiplication, then
$M_E(N)\ll_{E, \varepsilon}  N^\varepsilon$
for any $\varepsilon > 0$, and asks if the same might be true for general elliptic curves, but up to now,
no bound better that the trivial bound is known for curves without complex multiplication.

In~\cite{DS:2013}, we studied the average behavior of $M_E(N)$ taken over a family of elliptic curves.
More precisely, given integers $a,b$, we let $E_{a,b}$ denote the elliptic curve given by the
Weierstrass equation $y^2=x^3+ax+b$; we let $\C$ denote the multiset defined by
\begin{equation*}
\C=\C(A,B):=\{E_{a,b}: |a|\le A, |b|\le B, \Delta(E_{a,b})\ne 0\}.
\end{equation*}
Under a suitable conjecture (specifically Conjecture~\ref{short interval bdh conj} on
page~\pageref{short interval bdh conj} with any $\eta < 1/2$) we showed that there exists an absolutely bounded function $K(N)$
(see equation~\eqref{defn of K(N)} on page~\pageref{defn of K(N)}) such that if  $A,B$ are large enough with
respect to $N$, then
\begin{equation} \label{mainDS1}
\frac{1}{\#\C}\sum_{E\in\B}M_{E}(N)
\sim K(N)\frac{N}{\varphi(N)\log N}
\end{equation}
as $N \rightarrow \infty$.

Apart from the ``arithmetic factor" $K(N)N/\varphi(N)$, the above result agrees with a na\"ive probabilistic
model for $M_E(N)$ where one supposes that the values $\# E_p(\F_p)$ are uniformly distributed in the interval
$(N^-,N^+)$.  This is explained in detail in~\cite{DS:2013}.
The occurrence of the weight $\varphi(N)$ appearing in the denominator on the right hand side
of~\eqref{mainDS1} suggested to the authors that perhaps this is another example of phenomena which are
governed by the Cohen-Lenstra Heuristics~\cite{CL:1984,CL:1984-2}, which predict that random
 groups occur with probability inversely proportional to the size of their automorphism groups.
 The purpose of this paper is to explore this connection further by studying the function
\begin{equation*}
M_E(G):=\#\{p: E_p(\F_p)\isom G\},
\end{equation*}
where $G$ is a finite Abelian group.
Given an elliptic curve $E$, it is well known that
\begin{equation*}
E_p(\F_p)\isom \Z/N_1\Z\times\Z/N_1N_2\Z,
\end{equation*}
for some positive integers $N_1,N_2$ satisfying the Hasse bound: $|p+1-N_1^2N_2|<2\sqrt p$.
For much of this work, we will restrict to the case when $N_1$ and $N_2$ are both odd as this reduces the
number of special cases to consider.

As with our study of $M_E(N)$ in~\cite{DS:2013}, the restriction imposed by the Hasse bound means that any
prime counted by $M_E(G)$ must lie in a very short interval near $N=\#G=N_1^2N_2$.
In particular, all of the primes are of size $N$, lying in an interval of length $4 \sqrt{N}$.  Even the Riemann
Hypothesis does not guarantee the existence of a prime in such a short interval.
Thus our work here, as in~\cite{DS:2013}, requires a
conjecture (Conjecture~\ref{short interval bdh conj}) concerning the
distribution of primes of size $X$ in intervals of length $X^\eta$. The case $\eta=1$ corresponds to the
classical Barban-Davenport-Halberstam Theorem.
The precise statement of this conjecture can be found on page~\pageref{short interval bdh conj}
in Section~\ref{proofs of lemmas}.

Recall that the \textit{exponent} of a finite Abelian group is the size of its largest cyclic subgroup.
In particular, for groups of the form $G=\Z/N_1\Z\times\Z/N_1N_2\Z$, the exponent of $G$ is given by
$\exp(G)=N_1N_2$.
The following is our main result.

\begin{thm} \label{main}
Assume Conjecture \ref{short interval bdh conj} holds for some $\eta<1/2$.
Let $\alpha,\beta>0$ and fixed.  Then there exists a nonzero and absolutely bounded function $K(G)$
such that for every non-trivial, odd order group $G = \Z/N_1 \Z \times \Z/N_1 N_2 \Z$, we have that
\begin{equation*}
\frac{1}{\#\C}\sum_{E\in\C}M_E(G) =
	\left(K(G) + O \left( \frac{1}{(\log\# G)^{\beta}} \right)\right)\frac{\# G}{\# \Aut(G)\log(\#G)},
\end{equation*}
when $\# G \rightarrow \infty$
provided that $A,B\ge (\#G)^{1/2} (\log\#G)^{\beta+1}$, $AB \ge (\#G)^{3/2} (\log\#G)^{\beta+2}$ and
the exponent of $G$ satisfies $\frac{\#G}{(\log\#G)^\alpha}\le\exp(G)\le \#G$.
Furthermore, the function $K(G)$ is given as a product over primes $\ell$ by
\begin{equation}\label{defn of K(G)}
\begin{split}
K(G)&:=\prod_{\ell\dnd N}\left(1 - \frac{ \leg{N-1}{\ell}^2 \ell + 1}{(\ell-1)^2 (\ell+1)} \right)
	   \prod_{\ell\mid N}\left(1-\frac{1}{\ell(\ell-1)}\right)\\
&\quad\times\prod_{\ell\mid N_1}\left(1+\frac{1}{\ell(\ell^2-\ell-1)}\right)
             \prod_{\substack{\ell\mid N_1\\ \ell\dnd N_2}}\left(1+\frac{\leg{-N_2}{\ell}}{\ell(\ell-1)}\right),
\end{split}
\end{equation}
where $N=\#G=N_1^2N_2$ and $\leg{\cdot}{\ell}$ is the usual Kronecker symbol.
\end{thm}

Many results similar to Theorem~\ref{main} may be equivalently (and perhaps more naturally) stated as
results about counting isomorphism classes of elliptic curves over finite fields which possess some
desired property.  In this case, the desired property is having group of $\F_p$-points isomorphic to $G$.
We now restate Theorem~\ref{main} in such an equivalent form.

Let $M_p(G)$ denote the weighted number of isomorphism classes of elliptic curves defined over $\F_p$ with
group isomorphic to $G$.  That is,
\begin{equation}\label{defn of M_p(G)}
M_p(G):=\sum_{\substack{E/\F_p\\ E(\F_p)\isom G}}\frac{1}{\#\Aut(E)},
\end{equation}
where the sum is taken over all isomorphism classes of elliptic curves $E$ defined over $\F_p$ and
$\#\Aut(E)$ is the number of automorphisms of $E$ as a curve over $\F_p$.
It is important to distinguish between the similar notations $\Aut(E)$ and $\Aut(E(\F_p))$.
The former refers to the $\F_p$-automorphisms of $E$ as a curve, while the latter refers to the
automorphisms of $E(\F_p)$ as a group.
Now let
\begin{equation*}
M(G):=\sum_pM_p(G).
\end{equation*}
With this notation, Theorem~\ref{main} is equivalent to the following estimate for $M(G)$, the weighted number of
isomorphism classes of elliptic curves defined over any prime finite field with group of points isomorphic to $G$.
\begin{thm} \label{main-rephrased}
Assume that Conjecture~\ref{short interval bdh conj} holds for some $\eta<1/2$.
Let $\alpha,\beta>0$ and fixed.
Then for every non-trivial, odd order group $G=\Z/N_1\Z\times\Z/N_1N_2\Z$, we have
\begin{equation*}
M(G)=\left(K(G) + O \left( \frac{1}{(\log\# G)^{\beta}} \right)\right)\frac{(\# G)^2}{\# \Aut(G)\log(\#G)},
\end{equation*}
as $\# G \rightarrow \infty$
provided that the exponent of $G$ satisfies $\frac{\#G}{(\log\#G)^\alpha}\le\exp(G)\le \#G$.
\end{thm}
\begin{rmk}
 The proofs of Theorems~\ref{main} and~\ref{main-rephrased} do not really require that
 Conjecture~\ref{shortBDH} hold for a fixed $\eta<1/2$.  It is sufficient that it hold for intervals of length
 $Y=\sqrt X/(\log X)^{\beta+1}$.
\end{rmk}


The restriction that $\#G/(\log\#G)^\alpha\le\exp(G)\le \#G$ may seem a bit severe.
We believe that our results should hold in the range $(\#G)^{\frac{1}{2}+\epsilon}\le\exp(G)\le\#G$
for any fixed $\epsilon > 0$.  Proving this would require that we assume a conjecture similar to
Conjecture~\ref{short interval bdh conj} for primes in short intervals {\it and} in a fixed arithmetic progression.
Unconditionally, it is possible to obtain upper bounds of the correct order of magnitude in this larger range.
This is the subject of a forthcoming paper with V. Chandee and D. Koukoulopoulos in which we show that
\begin{equation*}
M(G)\ll\frac{(\# G)^2}{\# \Aut (G)\log(\#G)},
\end{equation*}
or equivalently,
\begin{equation*}
\frac{1}{\#\C}\sum_{E\in\C}M_E(G) \ll \frac{\# G}{\# \Aut (G)\log{(\# G)}},
\end{equation*}
both holding for $(\#G)^{\frac{1}{2}+\epsilon}\le\exp(G)\le\#G$.
However, lower bounds are impossible without hypothesis.

One should not expect our results to hold without some restriction on the size of of the exponent.
In particular, not all groups of the form $G= \Z/N_1 \Z \times \Z/N_1 N_2 \Z$ occur as the group of points on an
elliptic curve over a finite field.   For example, the authors of~\cite{BPS:2012} have noted that the group
$\Z/11\Z\times\Z/11\Z$ never occurs.
This is perhaps surprising at first since given a positive integer $N$ and a prime $p$ in the
range~\eqref{p range}, a theorem of Deuring~\cite{Deu:1941} ensures that there is always an elliptic curve
$E/\F_p$ possessing $N$ points.
Given a positive integer $N$, we believe that there should always be a prime close enough, but as we noted
earlier, this is not provable even under the Riemann Hypothesis.  A refinement of the Deuring result
(see Theorem~\ref{schoof thm} below) implies that given an order $N$ group
$G=\Z/N_1\Z\times\Z/N_1N_2\Z$ and a prime $p\equiv 1\pmod{N_1}$ in the range~\eqref{p range},
there is always an elliptic curve with $E(\F_p)\isom G$.
However, in the extreme case when $G = \Z/N_1 \Z \times \Z/N_1 \Z$, we are looking for a prime
$p\equiv 1\pmod{N_1}$ in the interval $(N_1^2 - 2N_1 + 1, N_1^2 + 2N_1 + 1)$, and we should not expect this
to happen very often as this interval contains exactly \textit{three} integers congruent to 1 modulo $N_1$.
In fact, letting $N_1$ and $N_2$ vary, it would seem very unlikely (though not impossible) to find an
elliptic curve $E/\F_p$ with $E(\F_p)\isom \Z/N_1\Z\times\Z/N_1N_2\Z$ unless $N_1$ grows slower than
an arbitrarily large (but fixed) power of $N_2$.  Note that this condition is equivalent to assuming
$\sqrt{N}/N_1\ge N^\epsilon$, which is equivalent to assuming $(\#G)^{\frac{1}{2}+\epsilon}\le\exp(G)$.

The remainder of the article is organized as follows.  We show in Section~\ref{reduce to class no avg} how
the proof of Theorem~\ref{main} is reduced to proving Theorem~\ref{main-rephrased}, and in turn, how the
proof of Theorem~\ref{main-rephrased} is reduced to the computation of a certain average of class numbers.
The computation of this average of class numbers occupies
Section~\ref{cond class no comp} (see Theorem~\ref{cond class number avg})
and Section~\ref{comp K_0} (see Proposition~\ref{compute euler prod}).
In Section~\ref{remove larger groups} we gather
together all of our intermediate results to complete the proof of Theorem~\ref{main-rephrased}.
Section~\ref{proofs of lemmas} contains the precise statement of
Conjecture~\ref{short interval bdh conj} as well as the proofs of a couple of auxiliary lemmas which are used
in Section~\ref{cond class no comp}.
Finally, in Section~\ref{conclusion} we close with some concluding remarks concerning the arithmetic factors $K(N)$ and $K(G)$,
which are respectively defined by~\eqref{defn of K(N)} and~\eqref{defn of K(G)}.

\section{Acknowledgement.}

The authors would like to thank Andrew Granville and Dimitris Koukoulopoulos for helpful discussions
while working on this project, and the anonymous referee for useful suggestions.
The first author also thanks Henri Cohen for enlightening discussions
about the Cohen-Lenstra Heuristics during the MSRI Program ``Arithmetic Statistics" in 2011.

The first author is supported by the Natural Sciences and Engineering Research Council
 of Canada [Discovery Grant 155635-2008] and the Fonds de recherche du Qu\'ebec - Nature et technologies [Octroi 166534].
 
 This work was completed while the second author was a postdoctoral fellow at Centre de recherches math\'ematiques in 
 Montr\'eal, Canada.  He is grateful for their support.

\section{Reduction to an average of class numbers.}\label{reduce to class no avg}

In this section we show how the proof of Theorem~\ref{main} is reduced to proving
Theorem~\ref{main-rephrased}.  We then show how the proof of Theorem~\ref{main-rephrased} is reduced to
computing a certain average of class numbers.



\begin{proof}[Proof that Theorem~\ref{main-rephrased} implies Theorem~\ref{main}]
For notational convenience, we let $N=\#G$.
The Hasse bound implies that any prime $p$ counted by $M_E(G)$ must fall in the range~\eqref{p range}.
Hence, interchanging the order of summation yields
\begin{equation}\label{interchange p and E}
\frac{1}{\#\C}\sum_{E\in\C}M_E(G)
=\frac{1}{\#\C}\sum_{N^-<p<N^+}\#\{E\in\C: E_p(\F_p)\isom G\}.
\end{equation}

To estimate the above, we group the $E$ in $\C$ according to which
$\F_p$-isomorphism class they reduce modulo $p$.  That is, we write
\begin{equation}\label{break over isom classes}
\#\{E\in\C: E_p(\F_p)\isom G\}
=\sum_{\substack{\tilde E/\F_p\\ \tilde E(\F_p)\isom G}}
	\#\{E\in\C: E_p(\F_p)\simeq_p\tilde E\},
\end{equation}
where the sum is over the $\F_p$-isomorphism classes of elliptic curves $\tilde E$ whose group of
$\F_p$-points is isomorphic to $G$.  We may assume that $N\ge 8$ so that
$N^-=(\sqrt N-1)^2>3$ and only primes greater than $3$ enter into the sum over $p$ above.
Therefore,  we may choose a model for each $\tilde E$ of the
form $y^2=x^3+\alpha x+\beta$, and a character sum argument as in~\cite[pp. 93-96]{FM:1996} yields
the estimate
\begin{equation}\label{reduction count}
\begin{split}
\#\{E\in\C: E_p(\F_p)\isom\tilde E\}
=
\frac{4AB}{\#\Aut(\tilde E)N}&+O\left(\frac{AB}{N^{2}}+\sqrt N(\log N)^2\right)\\
&+\begin{cases}
O\left(\frac{A\log N}{\sqrt N}+\frac{B\log N}{\sqrt N}\right)&\text{if }\alpha\beta\ne 0,\\
O\left(\frac{A\log N}{\sqrt N}+B\log N\right)&\text{if }\alpha=0,\\
O\left(A\log N+\frac{B\log N}{\sqrt N}\right)&\text{if }\beta=0
\end{cases}
\end{split}
\end{equation}
for $N^-<p<N^+$.

We now recall the definition of $M_p(G)$ as given by~\eqref{defn of M_p(G)}.
Then equations~\eqref{interchange p and E}, \eqref{break over isom classes}, and~\eqref{reduction count}
imply that
\begin{equation*}
\begin{split}
\frac{1}{\#\C}\sum_{E\in\C}M_E(G)
&=\left[
\frac{1}{N}
	+O\left(
	\frac{1}{N^{2}}
	+\frac{\log N}{\sqrt N}\left(\frac{1}{A}+\frac{1}{B}\right)
	+\frac{\sqrt N(\log N)^2}{AB}
	\right)
	\right]\sum_{N^-<p<N^+}M_p(G)\\
&\quad\quad+O\left(\left(\frac{1}{A}+\frac{1}{B}\right)\log N\sum_{N^-<p<N^+}1\right)
\end{split}
\end{equation*}
since $\#\C=4AB+O(A+B)$ and there are at most $10$ isomorphism classes of elliptic curves $\tilde E/\F_p$
with $\alpha\beta=0$.  Recalling that $M(G)=\sum_pM_p(G)$ and $\{N^-<p<N^+\}\ll\sqrt N/\log N$, we see
that Theorem~\ref{main} follows from Theorem~\ref{main-rephrased} provided that
$A,B\ge \sqrt{N}(\log N)^{\beta+1}$ and
$AB\ge N^{3/2}(\log N)^{\beta+2}$.
\end{proof}

Now let $N$ and $m$ be positive integers with $m^2$ dividing $N$, and define
\begin{equation*}
M_p(N;m):=\sum_{\substack{E/\F_p\\ \#E(\F_p)=N\\ E(\F_p)[m]\isom\Z/m\Z\times\Z/m\Z}}\frac{1}{\#\Aut(E)},
\end{equation*}
the weighted number of isomorphism classes of elliptic curves defined over $\F_p$ which have exactly $N$
points over $\F_p$ and whose $\F_p$-rational $m$-torsion subgroup is isomorphic to $\Z/m\Z\times\Z/m\Z$.
We also define
\begin{equation*}
M(N;m):=\sum_pM_p(N;m),
\end{equation*}
the weighted number of isomorphism classes of elliptic curves defined over any prime finite field which
have exactly $N$ points over $\F_p$ and whose $\F_p$-rational $m$-torsion subgroup is isomorphic to
$\Z/m\Z\times\Z/m\Z$.
We will compute $M(G)$ by first computing $M(N;m)$.
\begin{lma}\label{sieve lma}
Let $G=\Z/N_1\Z\times\Z/N_1N_2\Z$.  Then
\begin{equation*}
M(G)=\sum_{k^2\mid N_2}\mu(k)M(N_1^2N_2; kN_1),
\end{equation*}
where $\mu(k)$ denotes the usual M\"obius function.
\end{lma}
\begin{proof}
Inclusion-exclusion.
\end{proof}

We now explain how computing $M(N;m)$ is equivalent to computing a certain average of class numbers.
Given a negative discriminant $d$, we let $h(d)$ denote the
class number of the unique
imaginary quadratic order of discriminant $d$, and we let $w(d)$ denote the cardinality of its unit group.
The \textit{Kronecker class number} of discriminant $D$ is defined by
\begin{equation*}
H(D):=\sum_{\substack{f^2\mid D\\ \frac{D}{f^2}\equiv 0,1\pmod{4}}}\frac{h(D/f^2)}{w(D/f^2)}.
\end{equation*}

Given a positive integer $N$ and a prime $p$, we define the ``discriminant polynomial" $D_N(p)$ by
\begin{equation}\label{defn of D}
D_N(p):=(p+1-N)^2-4p.
\end{equation}
Adapting the proofs in~\cite[Lemma (4.8) and Theorem (4.9)]{Sch:1987}
to count isomorphism classes of elliptic curves weighted by the size of their automorphism groups,
we obtain the following.
\begin{thm}\label{schoof thm}
Let $p$ be a prime, $N$ a positive integer such that $|p+1-N|<2\sqrt p$, and $m$ a positive integer such that
$m^2\mid N$.  Then the weighted number of $\F_p$-isomorphism classes of elliptic curves having exactly $N$
points and $E(\F_p)[m]\isom\Z/m\Z\times\Z/m\Z$ is given by
\begin{equation*}
M_p(N;m)
=\begin{cases}
H\left(\frac{D_N(p)}{m^2}\right)&\text{if }m\mid p-1\text{ and } m^2\mid N,\\
0&\text{otherwise}.
\end{cases}
\end{equation*}
\end{thm}
\begin{rmk} It is easy to check that the conditions $m \mid p-1$ and $m^2 \mid N$ imply that $m^2 \mid D_N(p)$ in the above theorem.
\end{rmk}
\begin{rmk}
The corresponding results in~\cite{Sch:1987} are more general as Schoof does not restrict to finite fields of prime order.
\end{rmk}
As an immediate corollary, we have the following
\begin{cor}\label{reduction to class number avg}
Let $N$ and $m$ be positive integers with $m^2$ dividing $N$.  Then
\begin{equation*}
M(N;m)=\sum_{\substack{N^-<p<N^+\\ p\equiv 1\pmod m}}H\left(\frac{D_N(p)}{m^2}\right).
\end{equation*}
\end{cor}

\section{Conditional estimates for the average of class numbers.}\label{cond class no comp}

The following lemma will be useful for bounding various sums that appear in this section.
We postpone its proof until Section~\ref{proofs of lemmas}.

\begin{lma}\label{BT cor}
Suppose that $N$, $u$, and $v$ are positive integers with $u^2\mid N$.
Let $D_N(p)$ be as defined by equation~\eqref{defn of D}.
If $(X,X+Y]\subseteq (N^-,N^+)$, then uniformly for $uv<Y$,
\begin{equation*}
 \#\left\{X<p\le X+Y:\begin{array}{rl}p&\equiv 1\pmod u\\ D_N(p)&\equiv 0\pmod{u^2v}\\ (p,v)&=1\end{array}\right\}
\ll\frac{\sqrt v}{\varphi(uv)}\frac{Y}{\log(Y/uv)}.
\end{equation*}
\end{lma}

The main result of this section is the following conditional (under Conjecture~\ref{short interval bdh conj})
estimate for $M(N;m)$,
the weighted number of isomorphism classes of elliptic curves $E$ defined over any prime finite field which
have exactly $N$ points over $\F_p$ and such that $E(\F_p)$ contains a subgroup isomorphic to
$\Z/m\Z\times\Z/m\Z$.

\begin{thm}\label{cond class number avg} Let $\alpha, \beta > 0$ be fixed, and
assume that Conjecture~\ref{shortBDH} holds for $\eta<\frac{1}{2}$.
Then there exists a function $K_0(N,m)$ such that
for every pair of odd positive integers $N$ and $m$ with $m^2\mid N$,
\begin{equation*}
M(N;m)=\sum_{\substack{N^-<p<N^+\\ p\equiv 1\pmod m}}H\left(\frac{D_N(p)}{m^2}\right)
= K_0(N,m)\frac{N}{\log{N}}  + O \left( \frac{N}{\varphi(m^2)(\log N)^{\beta+1}} \right),
\end{equation*}
provided that $m\le (\log N)^\alpha$.  Furthermore, the function $K_0(N,m)$ is given by the absolutely
convergent sum
\begin{equation}\label{defn K0}
K_0(N,m):=\sum_{\substack{f =1\\ m\mid f\\ (f,2)=1}}^\infty\frac{1}{f}\sum_{n=1}^\infty
 	\frac{1}{n\varphi(4nf^2)}\sum_{\substack{a=1\\ a\equiv 1\pmod 4}}^{4n}\leg{a}{n}\#C_{N}(a,n,f)
\end{equation}
where
\begin{equation}\label{def big C}
C_{N}(a,n,f):=\left\{b\in(\Z/4nf^2\Z)^\times: D_N(b)\equiv af^2\pmod{4nf^2}\right\}.
\end{equation}
\end{thm}

\begin{proof}
We begin by using the definition of the Kronecker class number and the class number formula to write
\begin{equation*}
\sum_{\substack{N^-<p<N^+\\ p\equiv 1\pmod m}}H\left(\frac{D_N(p)}{m^2}\right)
=\frac{1}{2\pi}
	\sum_{\substack{N^-<p<N^+\\ p\equiv 1\pmod m}}
	\sum_{\substack{f^2\mid\frac{D_N(p)}{m^2}\\ d_{N,fm}(p)\equiv 0,1\pmod 4}}
	\frac{\sqrt{|D_N(p)|}}{fm}L\left(1,\chi_{d_{N,fm}(p)}\right),
\end{equation*}
where $\chi_d:=\leg{d}{\cdot}$ is the Kronecker symbol associated to the discriminant $d$,
 \begin{equation*}
 L(1,\chi_d)=\sum_{n\ge 1}\frac{\chi_d(n)}{n},
 \end{equation*}
and we write $d_{N,f}(p):=D_N(p)/f^2$ whenever $f^2\mid D_N(p)$.
We may assume that $N>5$.  As a result, the prime $2$ does not enter into the sum over $p$ above.
Since $N$ is odd, we have that $D_N(p)=(p+1-N)^2-4p\equiv 1\pmod 4$, and hence it follows that each
$f$ above is odd and $d_{N,fm}(p)\equiv 1\pmod 4$.  Therefore, we may omit the congruence condition under
the sum over $f$ above.   Furthermore, if $p\mid f$ and $f^2\mid D_N(p)$ it follows that $p=2$, but this is
contrary to our assumption that $N>5$.
Since $-3$ is the largest discriminant possible for an imaginary quadratic order, we have that
$d_{N,fm}(p)\le -3$, and it follows that if $(fm)^2\mid D_N(p)$, then $fm\le2\sqrt{N/3}$.  Therefore,
after rearranging the order of summation, we have that
\begin{equation}\label{unpack class numbers}
\sum_{\substack{N^-<p<N^+\\ p\equiv 1\pmod m}}H\left(\frac{D_N(p)}{m^2}\right)
=\frac{1}{2\pi}\sum_{\substack{f\le\frac{2}{m}\sqrt{\frac{N}{3}}\\ (f,2)=1}}\frac{1}{fm}
	\sum_{\substack{N^-<p<N^+\\ p\equiv 1\pmod m\\ (fm)^2\mid D_N(p)\\ p\nmid f}}
	\sqrt{|D_N(p)|}L\left(1,\chi_{d_{N,fm}(p)}\right).
\end{equation}

Let $V$ be a positive parameter to be chosen.  Using Lemma~\ref{BT cor} (with $u=m$ and $v=f$),
we see that the contribution made to the above by the values of $f$
which are larger than $V$ is bounded (up to a constant) by
\begin{equation*}
 \sqrt N\sum_{\substack{V<f\le\frac{2}{m}\sqrt{\frac{N}{3}}\\ (f,2)=1}}\frac{\log(4N/(fm)^2)}{fm}
	\sum_{\substack{N^-<p<N^+\\ p\equiv 1\pmod m\\ fm^2\mid D_N(p)\\ p\nmid f}}1
\ll N\sum_{f>V}\frac{1}{\sqrt f\varphi(fm^2)}
\ll\frac{N}{\varphi(m^2)\sqrt V}.
\end{equation*}
Choosing $V=(\log N)^{2(\beta+1)}$, we have that
\begin{equation}\label{trunc f sum}
\begin{split}
\sum_{\substack{N^-<p<N^+\\ p\equiv 1\pmod m}}H\left(\frac{D_N(p)}{m^2}\right)
&=\frac{1}{2\pi}\sum_{\substack{f\le V\\ (f,2)=1}}\frac{1}{fm}
	\sum_{\substack{N^-<p<N^+\\ p\equiv 1\pmod m\\ (fm)^2\mid D_N(p)\\ p\nmid f}}
	\sqrt{|D_N(p)|}L\left(1,\chi_{d_{N,fm}(p)}\right)\\
&\quad+O\left(\frac{N}{\varphi(m^2)(\log N)^{\beta+1}}\right).
\end{split}
\end{equation}
If $U>0$ and $d$ is a discriminant, then  Burgess' bound for character
sums~\cite[Theorem 2]{Bur:1963} and partial summation imply that
\begin{equation*}
L(1,\chi_d)=\sum_{n=1}^\infty\frac{\chi_d(n)}{n}
=\sum_{n\le U}\frac{\chi_d(n)}{n}
+O\left(\frac{|d|^{7/32}}{\sqrt U}\right).
\end{equation*}
Thus, if we truncate the $L$-series in~\eqref{trunc f sum} at $U$ and use the
Brun-Titchmarsh inequality~\cite[p.~167]{IK:2004},
we see that the contribution made by the tail of the $L$-series is bounded (up to a constant) by
\begin{equation*}
\frac{N^{23/32}}{\sqrt U}\sum_{f\le V}\frac{1}{(fm)^{1+7/16}}\sum_{\substack{N^-<p<N^+\\ p\equiv 1\pmod m}}1
\ll\frac{N^{39/32}}{\varphi(m^2)m^{7/16}\log(4\sqrt N/m)}
\end{equation*}
uniformly for $m\le\sqrt N$.
Choosing $U=N^{\frac{7}{16}}(\log N)^{2\beta}$, we have that
\begin{equation}\label{trunc n sum}
\begin{split}
\sum_{\substack{N^-<p<N^+\\ p\equiv 1\pmod m}}H\left(\frac{D_N(p)}{m^2}\right)
&=\frac{1}{2\pi}\sum_{\substack{f\le V\\ (f,2)=1}}\frac{1}{fm}
	\sum_{n\le U}\frac{1}{n}
	\sum_{\substack{N^-<p<N^+\\ p\equiv 1\pmod m\\ (fm)^2\mid D_N(p)\\ p\nmid f}}
	\leg{d_{N,fm}(p)}{n}\sqrt{|D_N(p)|}\\
&\quad+O\left(\frac{N}{\varphi(m^2)(\log N)^{\beta+1}}\right).
\end{split}
\end{equation}

At this point, it is convenient to break the interval $(N^-,N^+)$  into subintervals of length
$Y=\frac{\sqrt N}{\floor{(\log N)^\gamma}}$, where $\gamma$ is some fixed positive parameter to be chosen.
For each integer $k\in I:=[-2\sqrt N/Y,2\sqrt N/Y)\cap\Z$, we write
$X=X_k:=N+1+kY$ so that
\begin{eqnarray} \label{sumoverk}
\sum_{\substack{N^-<p<N^+\\ p\equiv 1\pmod m\\ (fm)^2\mid D_N(p)\\ p\nmid f}}
	\leg{d_{N,fm}(p)}{n}\sqrt{|D_N(p)|}
=\sum_{k\in I}
\sum_{\substack{X_k<p<X_k+Y\\ p\equiv 1\pmod m\\ (fm)^2\mid D_N(p)\\ p\nmid f}}
	\leg{d_{N,fm}(p)}{n}\sqrt{|D_N(p)|}.
\end{eqnarray}
On the interval $(X,X+Y]$, we approximate $\sqrt{|D_N(p)|}$ by $\frac{\sqrt{|D_N(X)|}\log p}{\log N}$.
Letting $X^*$ be the value of $t$ minimizing the function $\sqrt{|D_N(t)|}$ on the interval $[X,X+Y]$,
we find that
\begin{equation*}
\left|
\sqrt{|D_N(p)|}-\frac{\sqrt{|D_N(X)|}\log p}{\log N}
\right|
\ll\begin{cases}
N^{\frac{1}{4}}\sqrt Y+\frac{\sqrt{|D_N(X)|}}{\sqrt N\log N}&\text{if } N^{\pm}\in [X,X+Y],\\
\frac{Y\sqrt N}{\sqrt{|D_N(X^*)|}}+\frac{\sqrt{|D_N(X)|}}{\sqrt N\log N}&\text{otherwise}.
\end{cases}
\end{equation*}
Using Euler-Maclaurin summation as in~\cite{DS:2013}, it is easy to show
\begin{equation}\label{euler-maclaurin}
 \sum_{k \in I} \sqrt{|D_N(X_k)|} = \frac{2 \pi N}{Y} + O \left( \sqrt{N} \right)
\end{equation}
and
\begin{equation*}
\sum_{k \in I} \frac{1}{\sqrt{|D_N(X_k^*)|}} \ll  \frac{1}{Y}.
\end{equation*}
Employing these estimates and Lemma~\ref{BT cor} (with $u=m$ and $v=f^2$), we find that
\begin{equation*}
\begin{split}
\sum_{\substack{N^-<p<N^+\\ p\equiv 1\pmod m\\ (fm)^2\mid D_N(p)\\ p\nmid f}}
	\leg{d_{N,fm}(p)}{n}\sqrt{|D_N(p)|}
&=\frac{1}{\log N}\sum_{k\in I}\sqrt{|D_N(X_k)|}
\sum_{\substack{X_k<p\le X_k+Y\\ p\equiv 1\pmod m\\ (fm)^2\mid D_N(p)\\ p\nmid f}}\leg{d_{N,fm}(p)}{n}\log p\\
&\quad+O\left(\frac{Y\sqrt{N}}{\varphi(mf)\log(Y/mf^2)}\right)
\end{split}
\end{equation*}
provided that $mf<Y$.  Summing this error over $n$ and $f$, we find that the total error is
$O\left(\frac{Y\sqrt N\log U}{\varphi(m^2)(\log(Y/m))}\right)$ provided that $m<Y$.
Therefore, if $\gamma=\beta +1$, then
\begin{equation*}
\begin{split}
\sum_{\substack{N^-<p<N^+\\ p\equiv 1\pmod m}}H\left(\frac{D_N(p)}{m^2}\right)
&=\frac{1}{2\pi\log N}\sum_{k\in I}\sqrt{|D_N(X_k)|}\sum_{\substack{f\le V\\ (f,2)=1,\\ n\le U}}\frac{1}{fmn}
	\sum_{\substack{X_k<p\le X_k+Y\\ p\equiv 1\pmod m\\ (fm)^2\mid D_N(p)\\ p\nmid f}}
	\leg{d_{N,fm}(p)}{n}\log p\\
&\quad+O\left(\frac{N}{\varphi(m^2)(\log N)^{\beta+1}}\right)
\end{split}
\end{equation*}
for $m$ as large as $N^{\frac{1}{2}-\epsilon}$ (for any fixed $\epsilon>0$).
It is also convenient to remove those primes $p$ which divide $n$ from the innermost sum above.
Doing this introduces an error which is $O\left(m^{-1}\sqrt N(\log V)(\log U)^2\right)$.
Given our choice of $U$ and $V$, this can easily be absorbed into the error term for $m$
as large as $N^{\frac{1}{2}-\epsilon}$.

Now let
\begin{equation*}
S_k(N,m;U,V)
:=\sum_{\substack{f\le V\\ (f,2)=1,\\ n\le U}}\frac{1}{fmn}
	\sum_{\substack{X_k<p\le X_k+Y\\ p\equiv 1\pmod m\\ (fm)^2\mid D_N(p)\\ p\nmid fn}}
	\leg{d_{N,fm}(p)}{n}\log p
\end{equation*}
so that
\begin{equation*}
\sum_{\substack{N^-<p<N^+\\ p\equiv 1\pmod m}}H\left(\frac{D_N(p)}{m^2}\right)
=\frac{1}{2\pi\log N}\sum_{k\in I}\sqrt{|D_N(X_k)|}S_k(N,m;U,V)
	+O\left(\frac{N}{\varphi(m^2)(\log N)^{\beta+1}}\right).
\end{equation*}
Using the fact that the Kronecker symbol $\leg{\cdot}{n}$ is periodic modulo $4n$,
we may write
\begin{equation}
S_k(N,m;U,V)
=\sum_{\substack{f\le V\\ (f,2)=1,\\ n\le U}}\frac{1}{fmn}
	\sum_{\substack{a=1\\ a\equiv 1\pmod 4}}^{4n}\leg{a}{n}
	\sum_{\substack{X_k<p\le X_k+Y\\ p\equiv 1\pmod m\\ (fm)^2\mid D_N(p)\\ d_{N,fm}(p)\equiv a\pmod{4n}\\ p\nmid nf}}\log p.
\end{equation}
We now note that since $m^2\mid N$, the condition $(fm)^2\mid D_N(p)$ implies that $p\equiv 1\pmod m$.
Therefore, making the change of variables $mf\mapsto f$ and reorganizing the innermost sum over $p$, we obtain the identity
\begin{equation}
 S_k(N,m;U,V)
=\sum_{\substack{f\le mV\\ m\mid f\\ (f,2)=1,\\ n\le U}}\frac{1}{fn}
	\sum_{\substack{a=1\\ a\equiv 1\pmod 4}}^{4n}\leg{a}{n}
	\sum_{b\in C_{N}(a,n,f)}
	\sum_{\substack{X_k<p\le X_k+Y\\ p\equiv b\pmod{4nf^2}}}\log p,
\end{equation}
where
\begin{equation*}
C_{N}(a,n,f):=\left\{b\in(\Z/4nf^2\Z)^\times: D_N(b)\equiv af^2\pmod{4nf^2}\right\}.
\end{equation*}
We choose to approximate $S_k(N,m;U,V)$ by
\begin{equation*}
\tilde S_k(N,m;U,V):=\sum_{\substack{f\le mV\\ m\mid f\\ (f,2)=1,\\ n\le U}}\frac{1}{fn}
	\sum_{\substack{a=1\\ a\equiv 1\pmod 4}}^{4n}\leg{a}{n}\#C_{N}(a,n,f)\frac{Y}{\varphi(4nf^2)}.
\end{equation*}
It is at this point that we must impose the condition that $m\le (\log N)^\alpha$.
To bound the error in the above approximation, we use Cauchy-Schwarz and Conjecture~\ref{shortBDH} to
see that
\begin{equation*}
\begin{split}
\left|S-\tilde S\right|
&\le\sum_{\substack{f\le mV\\ m\mid f,\\ n\le U}}\frac{1}{fn}
	\sum_{b\in(\Z/4nf^2\Z)^\times}\left|\theta(X_k,Y;4nf^2,b)-\frac{Y}{\varphi(4nf^2)}\right|\\
&\le\sum_{\substack{f\le mV\\ m\mid f}}\frac{1}{f}
	\left[\sum_{n\le U}\frac{\varphi(4nf^2)}{n^2}\right]^{1/2}
	\left[\sum_{n\le U}\sum_{b\in(\Z/4nf^2\Z)^\times}\left|E(X_k,Y;4nf^2,b)\right|^2\right]^{1/2}\\
&\ll\frac{Y}{(\log N)^\upsilon}
\end{split}
\end{equation*}
for any choice of $\upsilon>0$ since
$4Um^2V^2\le 4N^{7/16}(\log N)^{2\alpha+6\beta+4}\ll Y/(\log N)^{\upsilon+1}$.
Summing this error over $k\in I$ and choosing $\upsilon=2\alpha+\beta$, we have that
\begin{equation*}
\begin{split}
 \sum_{\substack{N^-<p<N^+\\ p\equiv 1\pmod m}}H\left(\frac{D_N(p)}{m^2}\right)
&=\frac{1}{2\pi\log N}\sum_{k\in I}\sqrt{|D_N(X_k)|}\tilde S_k(N,m;U,V)\\
&\quad+O\left(\frac{N}{\varphi(m^2)(\log N)^{\beta+1}}\right)
\end{split}
\end{equation*}
since $m\le (\log N)^\alpha$.
Now, let
\begin{equation}
 K_0(N,m;U,V):=\sum_{\substack{f\le mV\\ m\mid f\\ (f,2)=1}}\frac{1}{f}\sum_{n\le U}
 	\frac{1}{n\varphi(4nf^2)}\sum_{\substack{a=1\\ a\equiv 1\pmod 4}}^{4n}\leg{a}{n}\#C_{N}(a,n,f).
\end{equation}
With this notation, we may rewrite the above as
\begin{equation}\label{ready to recomplete}
 \begin{split}
 \sum_{\substack{N^-<p<N^+\\ p\equiv 1\pmod m}}H\left(\frac{D_N(p)}{m^2}\right)
&=\frac{Y}{2\pi\log N}\sum_{k\in I}\sqrt{|D_N(X_k)|}K_0(N,m;U,V)\\
&\quad+O\left(\frac{N}{\varphi(m^2)(\log N)^{\beta+1}}\right).
\end{split}
\end{equation}
We now require the following lemma whose proof we postpone until Section~\ref{proofs of lemmas}.
\begin{lma}\label{trunc K0}
 For $\epsilon, U,V>0$,
\begin{equation*}
 K_0(N,m)=K_0(N,m;U,V)+O\left(\frac{N^\epsilon}{\varphi(m^2)\sqrt U}+\frac{\log\log N}{\varphi(m^2)V}\right).
\end{equation*}
\end{lma}
Given our choice of $U$ and $V$, this lemma together with~\eqref{euler-maclaurin}
and~\eqref{ready to recomplete} are sufficient to complete the proof of the theorem.
\end{proof}

%

\section{Computing the arithmetic factor $K_0(N,m)$.}\label{comp K_0}

The main result of this section is the following factorization of $K_0(N,m)$ as an
Euler product.

\begin{prop}\label{compute euler prod} Let $m$ and $N$ be odd positive integers with $m^2\mid N$, and
let $K_0(N,m)$ be as defined in Theorem~\ref{cond class number avg}. Then
\begin{equation*}
K_0(N,m) =\frac{N}{\varphi(N)m^2}K(N)K(N,m),
\end{equation*}
where
\begin{equation*}
K(N,m):=\prod_{\substack{\ell\mid m\\ 2\nmid\nu_\ell(N)}}
	 \left(\frac{\ell^{\nu_\ell(N)+1}-\ell^{2\nu_\ell(m)}}{\ell^{\nu_\ell(N)+1}-\ell^{\nu_\ell(N)}-1}\right)
	\prod_{\substack{\ell\mid m\\ 2\mid\nu_\ell(N)}}
	 \left(\frac{\ell^{\nu_\ell(N)+2}-\ell^{2\nu_\ell(m)+1}+\leg{-N_{(\ell)}}{\ell}\ell^{2\nu_\ell(m)}}
	{\ell^{\nu_\ell(N)+2}-\ell^{\nu_\ell(N)+1}-\ell+\leg{-N_{(\ell)}}{\ell}}\right)
\end{equation*}
and
\begin{equation}\label{defn of K(N)}
K(N)\!:=\!\prod_{\ell\nmid N}
		\!\left(1-\frac{\leg{N-1}{\ell}^2\ell+1}{(\ell+1)(\ell-1)^2}\right)
	\!\!\!\prod_{\substack{\ell\mid N\\ 2\nmid\nu_\ell(N)}}
		\!\!\!\!\left(1-\frac{1}{\ell^{\nu_\ell(N)}(\ell-1)}\right)
	\!\!\!\prod_{\substack{\ell\mid N\\ 2\mid\nu_\ell(N)}}
		\!\!\!\!\left(1-\frac{\ell-\leg{-N_{(\ell)}}{\ell}}{\ell^{\nu_\ell(N)+1}(\ell-1)}\right)\!.
\end{equation}
Here $\nu_\ell(N)$ is the usual $\ell$-adic valuation, $N_{(\ell)}:=N/\ell^{\nu_\ell(N)}$ denotes the $\ell$-free part
of $N$, and $\leg{\cdot}{\ell}$ is the usual Kronecker symbol.
\end{prop}

\begin{rmk}
It will be convenient in the next section to note that $K(N)$ is absolutely bounded as a function of $N$,
$N/\varphi(N)\ll\log\log N$, and $K(N,m) \ll m/\varphi(m) \ll\log\log m$.
\end{rmk}

We note that $K_0(N)=K_0(N,1)$ was computed in~\cite{DS:2013}.
We will appeal often to results from~\cite{DS:2013} in our computation of $K_0(N,m)$.

\begin{proof}[Proof of Proposition~\ref{compute euler prod}]
We begin by using the Chinese Remainder Theorem to write
\begin{equation*}
K_0(N,m):=
	\sum_{\substack{f=1\\ m\mid f\\ (f,2)=1}}^\infty\frac{1}{f}
	\sum_{n=1}^\infty\frac{1}{n\varphi(4nf^2)}
	\sum_{\substack{a\in(\Z/4n\Z)^*\\ a\equiv 1\pmod 4}}\leg{a}{n}
	\prod_{\ell\mid 4nf}\#C_N^{(\ell)}(a,n,f),
\end{equation*}
where for each prime $\ell$ dividing $4nf$ we write
\begin{equation*}
C_N^{(\ell)}(a,n,f):=
\left\{z\in(\Z/\ell^{\nu_\ell(4nf^2)}\Z)^*:
	D_N(z)\equiv af^2\pmod{\ell^{\nu_\ell(4nf^2)}}\right\}.
\end{equation*}
We note that if $\ell$ is a prime dividing $f$ but not $n$, then
$af^2\equiv 0\pmod{\ell^{\nu_\ell(4nf^2)}}$.  Hence, it is clear that
$\#C_N^{(\ell)}(a,n,f)=\#C_N^{(\ell)}(1,1,f)$ if $\ell\mid f$ but $\ell\dnd n$.
In~\cite[Lemma 10]{DS:2013} it was shown that
\begin{equation*}
\#C_N^{(2)}(a,n,f)=2\mathcal S_2(n,a),
\end{equation*}
where
\begin{equation*}
\mathcal S_2(n,a)=
\begin{cases}
2&\text{if }\nu_2(4nf^2)=2+\nu_2(n)=2,\\
4&\text{if }\nu_2(4nf^2)=2+\nu_2(n)\ge 3\text{ and } a\equiv 5\pmod 8,\\
0&\text{otherwise}.
\end{cases}
\end{equation*}
Therefore, letting $n'$ denote the odd part of $n$ and
\begin{equation*}
c_{N,f}(n):=\sum_{\substack{a\in(\Z/4n\Z)^*\\ a\equiv 1\pmod 4}}
\leg{a}{n}\mathcal S_2(n,a)
\prod_{\ell\mid n'}\#C_N^{(\ell)}(a,n,f),
\end{equation*}
we may write
\begin{equation}\label{sep n and f}
\begin{split}
K_0(N,m)
&=\sum_{\substack{f=1\\ m\mid f\\ (f,2)=1}}^\infty\frac{1}{f^2}
	\sum_{n=1}^\infty\frac{2}{n\varphi(4nf)}
	\left[\prod_{\substack{\ell\mid f\\ \ell\dnd n}}\#C_N^{(\ell)}(a,n,f)\right]c_{N,f}(n)\\
&=\sum_{\substack{f=1\\ m\mid f\\ (f,2)=1}}^\infty\strut^\prime\frac{\prod_{\ell\mid f}
	\#C_N^{(\ell)}(1,1,f)}{f^2\varphi(f)}
	\sum_{n=1}^\infty\frac{2\varphi((n,f))}{(n,f)n\varphi(4n)}
	\left[\prod_{\ell\mid (f,n)}\#C_N^{(\ell)}(1,1,f)\right]^{-1}c_{N,f}(n),
\end{split}
\end{equation}
where the $\prime$ on the sum over $f$ is meant to indicate that the sum is to be restricted
to those $f$ that are not divisible by any prime $\ell$ for which $\#C_N^{(\ell)}(1,1,f)=0$.
The following was shown in~\cite{DS:2013}, and we state it here without proof.
\begin{lma}\label{compute factors of c}
Suppose that $N$ and $f$ are odd.
The function $c_{N,f}(n)$ is multiplicative in $n$.
Let $\alpha$ be a positive integer and $\ell$ an odd prime.
Then
\begin{equation*}
\frac{c_{N,f}(2^\alpha)}{2^{\alpha-1}}=(-1)^\alpha 2.
\end{equation*}
If $\ell\mid f$ and $\ell\dnd N$, then
\begin{equation*}
\frac{c_{N,f}(\ell^\alpha)}{\ell^{\alpha-1}}
=\#C_N^{(\ell)}(1,1,f)
\begin{cases}
\ell-1&\text{if }2\mid\alpha,\\
0&\text{if }2\dnd\alpha.
\end{cases}
\end{equation*}
If $\ell\mid N$ and $\ell\dnd f$, then
\begin{equation*}
\frac{c_{N,f}(\ell^\alpha)}{\ell^{\alpha-1}}=\ell-2.
\end{equation*}
If $\ell\dnd Nf$, then
\begin{equation*}
\begin{split}
\frac{c_{N,f}(\ell^\alpha)}{\ell^{\alpha-1}}
&=\begin{cases}
\ell-1-\leg{N}{\ell}-\leg{N-1}{\ell}^2
	&\text{if }2\mid\alpha,\\
-1-\leg{N-1}{\ell}^2
	&\text{if }2\dnd\alpha.
\end{cases}
\end{split}
\end{equation*}
If $\ell\mid (f,N)$ and $2\nu_\ell(f)<\nu_\ell(N)$, then
\begin{equation*}
\frac{c_{N,f}(\ell^\alpha)}{\ell^{\alpha-1}}
=\#C_N^{(\ell)}(1,1,f)(\ell-1).
\end{equation*}
If $\ell\mid (f,N)$ and $\nu_\ell(N)<2\nu_\ell(f)$, then
\begin{equation*}
\frac{c_{N,f}(\ell^\alpha)}{\ell^{\alpha-1}}
=\#C_N^{(\ell)}(1,1,f)\begin{cases}
\ell-1&\text{if }2\mid \alpha,\\
0&\text{if }2\dnd\alpha.
\end{cases}
\end{equation*}
If $\ell\mid (f,N)$ and $\nu_\ell(N)=2\nu_\ell(f)$, then
\begin{equation*}
\frac{c_{N,f}(\ell^\alpha)}{\ell^{\alpha-1}}
=\#C_N^{(\ell)}(1,1,f)\begin{cases}
\left(\ell-1-\leg{N_{(\ell)}}{\ell}+\leg{-N_{(\ell)}}{\ell}\right)&\text{if }2\mid\alpha,\\
\left(\leg{-N_{(\ell)}}{\ell}-1\right)&\text{if }2\dnd\alpha,\\
\end{cases}
\end{equation*}
where $N_{(\ell)}=N/\ell^{\nu_\ell(N)}$ denotes the $\ell$-free part of $N$.
\end{lma}

Using this result, the sum over $n$ in equation~\eqref{sep n and f} may be factored as
\begin{equation*}
\begin{split}
\sum_{n=1}^\infty\frac{2\varphi((n,f))}{(n,f)n\varphi(4n)}
	\left[\prod_{\ell\mid (f,n)}\#C_N^{(\ell)}(1,1,f)\right]^{-1}c_{N,f}(n)
 =\frac{2}{3}
	\prod_{\substack{\ell\dnd f\\ \ell\mid N}}F_0(\ell)	
	\prod_{\ell\dnd 2fN}F_1(\ell)
	\prod_{\ell\mid f}F_2(\ell,f),
\end{split}
\end{equation*}
where for any odd prime $\ell$, we make the definitions
\begin{align*}
F_0(\ell)&:=\left(1+\frac{\ell-2}{(\ell-1)^2}\right);\\
F_1(\ell)
&:=\left(1-
	\frac{\leg{N-1}{\ell}^2\ell
	+\leg{N}{\ell}+\leg{N-1}{\ell}^2+1}{(\ell-1)(\ell^2-1)}
\right);\\
F_2(\ell,f)&:=\begin{cases}
\displaystyle
\left(1+\frac{1}{\ell(\ell+1)}\right)&\text{if }\nu_\ell(N)<2\nu_\ell(f),\\
\displaystyle
\left(1+\frac{1}{\ell}\right)&\text{if }\nu_\ell(N)>2\nu_\ell(f),\\
\displaystyle
\left(1+
	 \frac{\leg{-N_{(\ell)}}{\ell}\ell+\leg{-N_{(\ell)}}{\ell}-\leg{N_{(\ell)}}{\ell}-1}{\ell(\ell^2-1)}
	\right)
	&\text{if }\nu_\ell(N)=2\nu_\ell(f).
\end{cases}
\end{align*}
Substituting this back into equation~\eqref{sep n and f} and rearranging slightly, we have that
\begin{equation}\label{f sum remaining}
K_0(N,m)
=\frac{2}{3}
	\prod_{\ell\mid N}F_0(\ell)
	\prod_{\ell \dnd 2N} F_1(\ell)
	\sum_{\substack{f=1\\ m\mid f\\ (f,2)=1}}^{\infty}\strut^\prime
		\frac{\prod_{\ell\mid f}\#C_N^{(\ell)}(1,1,f)}{\varphi(f)f^2}
		\prod_{\ell\mid (f,N)}\frac{F_2(\ell,f)}{F_0(\ell)}
		\prod_{\substack{\ell\mid f\\ \ell\dnd N}}\frac{F_2(\ell,f)}{F_1(\ell)}.
\end{equation}
In~\cite{DS:2013}, we showed that for any odd prime power $\ell^\alpha$,
\begin{eqnarray} \label{valueofnumberC}
\#C_N^{(\ell)}(1,1,\ell^\alpha)=\begin{cases}
1+\leg{N(N-1)^2}{\ell}&\text{if }\ell\dnd N,\\
2\ell^{\nu_\ell(N)/2}&\text{if }1\le\nu_\ell(N)<2\alpha,\ 2\mid\nu_\ell(N),
	\text{ and }\leg{N_{(\ell)}}{\ell}=1,\\
\ell^{\alpha}&\text{if }2\alpha\le\nu_\ell(N),\\
0&\text{otherwise}.
\end{cases}
\end{eqnarray}
Hence, the sum over $f$ in equation~\eqref{f sum remaining} may be factored as
\begin{equation*}
\begin{split}
&\sum_{\substack{f=1\\ m\mid f\\ (f,2)=1}}^{\infty}\strut^\prime
		\frac{\prod_{\ell\mid f}\#C_N^{(\ell)}(1,1,f)}{\varphi(f)f^2}
		\prod_{\ell\mid (f,N)}\frac{F_2(\ell,f)}{F_0(\ell)}
		\prod_{\substack{\ell\mid f\\ \ell\dnd N}}\frac{F_2(\ell,f)}{F_1(\ell)}
=\prod_{\ell\nmid 2N}F_3(\ell)
\prod_{\substack{\ell\mid N\\ \ell\nmid m}}F_4(\ell)
\prod_{\ell\mid m}F_5(\ell),
\end{split}
\end{equation*}
where for any odd prime $\ell$, we make the definitions
\begin{align*}
F_3(\ell)&:=1+\frac{1+\leg{N(N-1)^2}{\ell}}{F_1(\ell)(\ell+1)(\ell-1)^2};\\
F_4(\ell)&:=\begin{cases}\displaystyle
1+\frac{\ell^{\nu_\ell(N)}-\ell}{F_0(\ell)\ell^{\nu_\ell(N)}(\ell-1)^2}&\text{if }2\nmid \nu_\ell(N),\\
\displaystyle
1+\frac{\ell^{\nu_\ell(N)}-\ell+\leg{-N_{(\ell)}}{\ell}}{F_0(\ell)\ell^{\nu_\ell(N)}(\ell-1)^2}&\text{if }2\mid\nu_\ell(N);
\end{cases}\\
F_5(\ell)&:=\frac{1}{F_0(\ell)\ell^{2\nu_\ell(m)}}\begin{cases}
\displaystyle
\frac{\ell(\ell^{\nu_\ell(N)}-\ell^{2\nu_\ell(m)})}{\ell^{\nu_\ell(N)}(\ell-1)^2}
	&\text{if }2\nmid\nu_\ell(N),\\
\displaystyle
\frac{\ell^{\nu_\ell(N)+2}-\ell^{2\nu_\ell(m)+1}+\leg{-N_{(\ell)}}{\ell}\ell^{2\nu_\ell(m)}}{\ell^{\nu_\ell(N)}(\ell-1)^2}
	&\text{if }2\mid \nu_\ell(N)
.
\end{cases}
\end{align*}
Substituting this back into equation~\eqref{f sum remaining}, we have that
\begin{equation}\label{ready for final rearrange}
\begin{split}
K_0(N,m)
&=\frac{2}{3}\prod_{\ell\nmid 2N}F_1(\ell)F_3(\ell)
	\prod_{\substack{\ell\mid N\\ \ell\nmid m}}F_0(\ell)F_4(\ell)
	\prod_{\ell\mid m}F_0(\ell)F_5(\ell)\\
&=\frac{2}{3}\prod_{\ell\nmid 2N}F_1(\ell)F_3(\ell)
	\prod_{\ell\mid N}F_0(\ell)F_4(\ell)
	\prod_{\ell\mid m}\frac{F_5(\ell)}{F_4(\ell)}.
\end{split}
\end{equation}
The result now follows by simplifying the factors.
\end{proof}

\section{Removing the larger group structures.}\label{remove larger groups}

In this section we complete the proof of Theorem~\ref{main-rephrased}.
Before doing so, we require the following lemma.
\begin{lma}\label{size of Aut lma}
Let $G=\Z/N_1\Z\times\Z/N_1N_2\Z$, and write $N=\#G=N_1^2N_2$.  Then
\begin{equation*}
\frac{\#G}{\#\Aut(G)}
=\frac{N}{\varphi(N)N_1^2}
	\prod_{\substack{\ell\mid N_1\\ \ell\nmid N_2}}\frac{\ell^2}{\ell^2-1}
	\prod_{\ell\mid (N_1,N_2)}\frac{\ell}{\ell-1}
\end{equation*}
\end{lma}
\begin{proof}
See~\cite[Lemma 2.1 and Theorem 4.1]{HR:2007}.
\end{proof}
\begin{rmk}
From the above factorization it is easy to see that $\#G/\#\Aut(G)> 1/\varphi(N_1^2)$.
\end{rmk}

We now combine all of our intermediate results to give the proof of Theorem~\ref{main-rephrased}.

\begin{proof}[Proof of Theorem~\ref{main-rephrased}]
Let $N=N_1^2N_2=\#G$.
Then the condition $\exp(G)\ge\#G/(\log\#G)^\alpha$ is equivalent to
$N_1\le(\log N)^\alpha$.
Now choose any $\gamma>\beta+1$.
By Lemma~\ref{sieve lma}, we have
\begin{equation*}
M(G)=\sum_{\substack{k^2\mid N_2\\ k\le (\log N)^{\gamma}}}\mu(k)M(N;kN_1)
	+\sum_{\substack{k^2\mid N_2\\ k> (\log N)^{\gamma}}}\mu(k)M(N;kN_1).
\end{equation*}
In~\cite[p.~656]{Len:1987}, we find the standard bound $H(D)\ll\sqrt{|D|}\log|D|(\log\log|D|)^2$.
Using this, Corollary~\ref{reduction to class number avg}, and the
Brun-Titchmarsh inequality~\cite[p.~167]{IK:2004}, we have that
\begin{equation*}
\begin{split}
\sum_{\substack{k^2\mid N_2\\ k> (\log N)^{\gamma}}}\mu(k)M(N;kN_1)
&\ll\sum_{\substack{k^2\mid N_2\\ k> (\log N)^{\gamma}}}
	\sum_{\substack{N^-<p<N^+\\ p\equiv 1\pmod{kN_1}}}\frac{\sqrt{N}\log N(\log\log N)^2}{kN_1}\\
&\ll\frac{N\log N(\log\log N)^2}{\varphi(N_1^2)}\sum_{k>(\log N)^\gamma}\frac{1}{\varphi(k^2)\log(\sqrt N/kN_1)}\\
&\ll\frac{N}{\varphi(N_1^2)(\log N)^{\beta+1}}.
\end{split}
\end{equation*}
Applying Theorem~\ref{cond class number avg} for the small $k$, we find that
\begin{equation*}
M(G)=\frac{N}{\log N}\sum_{\substack{k^2\mid N_2\\ k\le (\log N)^{\gamma}}}\mu(k)K_0(N,kN_1)
	+O\left(\frac{N}{\varphi(N_1^2)(\log N)^{\beta+1}}\right).
\end{equation*}
By the remark following the statement of Proposition~\ref{compute euler prod},
\begin{equation*}
\sum_{\substack{k^2\mid N_2\\ k>(\log N)^\gamma}}\mu(k)K_0(N,kN_1)
\ll\sum_{k>(\log N)^{\gamma}}\frac{(\log\log N)^2}{k^2N_1^2}
\ll\frac{1}{\varphi(N_1^2)(\log N)^{\beta+1}},
\end{equation*}
and therefore,
\begin{equation}\label{ready to compute final arith factor}
M(G)=\frac{N}{\log N}\sum_{k^2\mid N_2}\mu(k)K_0(N,kN_1)
	+O\left(\frac{N}{\varphi(N_1^2)(\log N)^{\beta+1}}\right).
\end{equation}
Using Proposition~\ref{compute euler prod}, we have that
\begin{equation}\label{prep k sum}
\sum_{k^2\mid N_2}\mu(k)K_0(N,kN_1)
=\frac{N}{\varphi(N)N_1^2}K(N)\sum_{k^2\mid N_2}\frac{\mu(k)}{k^2}K(N,kN_1).
\end{equation}
Now let $K^{(\ell)}(N,m)$ stand for the factor of $K(N,m)$ coming from the prime $\ell$, i.e.,
\begin{equation*}
K^{(\ell)}(N,m):=
\begin{cases}
\displaystyle
\frac{\ell^{\nu_\ell(N)+1}-\ell^{2\nu_\ell(m)}}{\ell^{\nu_\ell(N)+1}-\ell^{\nu_\ell(N)}-1}
	&\text{if }\ell\mid m\text{ and }2\nmid\nu_\ell(N),\\
\displaystyle
\frac{\ell^{\nu_\ell(N)+2}-\ell^{2\nu_\ell(m)+1}+\leg{-N_{(\ell)}}{\ell}\ell^{2\nu_\ell(m)}}
	{\ell^{\nu_\ell(N)+2}-\ell^{\nu_\ell(N)+1}-\ell+\leg{-N_{(\ell)}}{\ell}}
	&\text{if }\ell\mid m\text{ and } 2\mid\nu_\ell(N).
\end{cases}
\end{equation*}
Then by multiplicativity, we have that
\begin{equation*}
\begin{split}
\sum_{k^2\mid N_2}\frac{\mu(k)}{k^2}K(N,kN_1)
&=\prod_{\ell\mid N_1}K^{(\ell)}(N,N_1)
	\sum_{k^2\mid N_2}\frac{\mu(k)}{k^2}
	\frac{\prod_{\ell\mid k}K^{(\ell)}(N,kN_1)}{\prod_{\ell\mid (k,N_1)}K^{(\ell)}(N,N_1)}\\
&=\prod_{\ell\mid N_1}K^{(\ell)}(N,N_1)
	\prod_{\substack{\ell^2\mid N_2\\ \ell\nmid N_1}}\left(1-\frac{K^{(\ell)}(N,\ell N_1)}{\ell^2}\right)
	\prod_{\substack{\ell^2\mid N_2\\ \ell\mid N_1}}
		\left(1-\frac{K^{(\ell)}(N,\ell N_1)}{\ell^2K^{(\ell)}(N,N_1)}\right)\\
&=\prod_{\substack{\ell\mid N_1\\ \ell^2\nmid N_2}}K^{(\ell)}(N,N_1)
	\prod_{\substack{\ell^2\mid N_2\\ \ell\nmid N_1}}\left(1-\frac{K^{(\ell)}(N,\ell)}{\ell^2}\right)
	\prod_{\substack{\ell^2\mid N_2\\ \ell\mid N_1}}\left(K^{(\ell)}(N,N_1)-\frac{K^{(\ell)}(N,\ell N_1)}{\ell^2}\right).
\end{split}
\end{equation*}
Recalling the definition of $K(N)$ as given by equation~\eqref{defn of K(N)}, we find that
\begin{equation*}
\begin{split}
K(N)\sum_{k^2\mid N_2}\frac{\mu(k)}{k^2}K(N,kN_1)
&=\prod_{\ell\nmid N}\left(1-\frac{\leg{N-1}{\ell}^2\ell+1}{(\ell+1)(\ell-2)^2}\right)
	\prod_{\substack{\ell\mid N_1\\ \ell\nmid N_2}}\frac{\ell^2}{\ell^2-1}F_6(\ell)
	\prod_{\substack{\ell\mid N_1\\ \ell\mid N_2}}\frac{\ell}{\ell-1}F_7(\ell)
	\prod_{\substack{\ell\nmid N_1\\ \ell\mid N_2}}F_8(\ell),
\end{split}
\end{equation*}
where for each odd prime $\ell$, we make the definitions
\begin{align*}
F_6(\ell)&:=1+\frac{(\leg{-N_{(\ell)}}{\ell}-1)\ell+\leg{-N_{(\ell)}}{\ell}}{\ell^3},
\\
F_7(\ell)&:=1-\frac{1}{\ell^2},
\\
F_8(\ell)&:=1-\frac{1}{\ell(\ell-1)}.
\end{align*}
Substituting the above into equation~\eqref{prep k sum} and using Lemma~\ref{size of Aut lma}, we find
(after some slight rearrangement of the factors) that
\begin{equation*}
\sum_{k^2\mid N_2}\mu(k)K_0(N,kN_1)=K(G)\frac{\#G}{\#\Aut(G)},
\end{equation*}
where $K(G)$ is defined by equation~\eqref{defn of K(G)}.
The result now follows by substituting this into equation~\eqref{ready to compute final arith factor} and using
the remark following the statement of Lemma~\ref{size of Aut lma}.
\end{proof}

\section{Conjecture \ref{short interval bdh conj} and proofs of Lemmas \ref{BT cor}
and \ref{trunc K0}.}\label{proofs of lemmas}

In this section, we give the precise statement of the conjecture for the distribution of primes
in short intervals that is needed to prove Theorems~\ref{main} and~\ref{main-rephrased}.
We also present the proofs of Lemmas~\ref{BT cor} and~\ref{trunc K0}.

Given real parameters $X, Y > 0$ and integers $q$ and $a$, we let $\theta(X,Y ; q,a)$ denote the
logarithmic weighted prime counting function
$$\theta(X, Y ; q, a) :=
\sum_{{X<p<X+Y} \atop {p \equiv a \pmod{q}}} \log{p},$$
and we let $E(X, Y ; q, a)$ be the error in approximating $\theta(X,Y ; q, a)$ by $Y/\varphi(q)$. That is,
$$E(X, Y ; q, a) := \theta(X, Y ; q, a) -Y/\varphi(q).$$

\begin{conj}\label{short interval bdh conj}
(Barban-Davenport-Halberstam for intervals of length $X^\eta$)
\label{shortBDH} Let $0 < \eta \leq 1$, and let $\beta > 0$ be arbitrary. Suppose
that $X^\eta \leq Y \leq X$, and that $Y / (\log{X})^\beta \leq Q \leq Y.$
Then
\begin{align*}
\sum_{q \leq Q} \sum_{{a=1}\atop{(a,q)=1}}^{q} \left| E(X,Y;q,a) \right|^2
\ll YQ \log{X}.
\end{align*}
\end{conj}
\begin{rmk}
If $\eta=1$, this is essentially the classical Barban-Davenport-Halberstam Theorem.
See for example~\cite[p.~196]{Dav:1980}.
The best results known are due to Languasco, Perelli, and Zaccagnini~\cite{LPZ:2010}, who
show that Conjecture~\ref{shortBDH} holds unconditionally for any
$\eta > 7/12$ and for any $\eta>1/2$ under the Generalized Riemann Hypothesis.
For the proofs of Theorems~\ref{main} and~\ref{main-rephrased}, we essentially need to assume that it
holds for some $\eta<1/2$.
\end{rmk}

\begin{proof}[Proof of Lemma~\ref{BT cor}]
First, we make the definition
\begin{equation}\label{def Delta}
 \Delta_{N,u}(l):=(l-N/u)^2-4N/u^2,
\end{equation}
and note that this is a quadratic polynomial in $l$ with integer coefficients since $u^2\mid N$.
We also note that if $p=1+lu$, then $\Delta_{N,u}(l)=D_N(p)/u^2$.
Therefore,
\begin{equation*}
\begin{split}
&\#\left\{X<p\le X+Y:\begin{array}{rl}p&\equiv 1\pmod u\\ D_N(p)&\equiv 0\pmod{u^2v}\\ (p,v)&=1\end{array}\right\}\\
&\quad=\#\left\{X<p\le X+Y:
	\begin{array}{rl} p&\equiv 1\pmod u\\
	\Delta_{N,u}\left(\frac{p-1}{u}\right)&\equiv 0\pmod{v}\\
	(p,v)&=1\end{array}\right\}\\
&\quad=\sum_{\substack{1\le l\le v\\ (1+lu,v)=1\\ \Delta_{N,u}(l)\equiv 0\pmod v}}
	\#\left\{X<p<X+Y: p\equiv 1+lu\pmod{uv}\right\}.
\end{split}
\end{equation*}
Using the Chinese Remainder Theorem, it easy (though perhaps a bit tedious) to show that
\begin{equation*}
\#\{l\in\Z/v\Z: \Delta_{N,u}(l)\equiv 0\pmod v\}\le 8\sqrt v.
\end{equation*}
We refer the reader to~\cite[Lemma 12]{DS:2013} where the same bound is shown for the polynomial
$D_N(p)$.  However, the same proof goes through for any monic quadratic with integer coefficients.
The above inequality together with the Brun-Titchmarsh inequality~\cite[p.~167]{IK:2004} implies that
\begin{equation}\label{bound p sum}
\#\left\{X<p\le X+Y:
	\begin{array}{rl} p&\equiv 1\pmod u\\ D_N(p)&\equiv 0\pmod{u^2v}\\ (p,v)&=1\end{array}\right\}
\ll\frac{\sqrt v}{\varphi(uv)}\frac{Y}{\log(Y/uv)},
\end{equation}
uniformly for $uv\le Y$.
\end{proof}

\begin{proof}[Proof of Lemma~\ref{trunc K0}]
In~\cite{DS:2013}, we showed that
\begin{equation*}
c_{N,f}(n)\ll\frac{n\prod_{\ell\mid (f,n)}\#C_N(1,1,f)}{\kappa_{2N}(n)},
\end{equation*}
where for any integer $m$, $\kappa_m(n)$ is the multiplicative function defined on prime
powers by
\begin{equation}\label{defn of kappa}
\kappa_m(\ell^\alpha)
:=\begin{cases}
\ell&\text{if }2\dnd\alpha\text{ and }\ell\dnd m,\\
1&\text{otherwise}.
\end{cases}
\end{equation}
Therefore,
\begin{equation}\label{diff bn K and truncated K}
\begin{split}
& K_0(N,m)-K_0(N,m;U,V) \\
&\quad\quad\ll\sum_{\substack{f\le mV\\ m\mid f}}\strut^\prime\frac{\prod_{\ell\mid f}\#C_N^{(\ell)}(1,1,f)}{f^2\varphi(f)}
	\sum_{n>U}\frac{2\varphi((n,f))c_{N,f}(n)}{(n,f)n\varphi(4n)
	\prod_{\ell\mid (n,f)}\#C_N^{(\ell)}(1,1,f)}\\
	&\quad\quad\quad+\sum_{\substack{f> mV\\ m\mid f}}\strut^\prime
		\frac{\prod_{\ell\mid f}\#C_N^{(\ell)}(1,1,f)}{f^2\varphi(f)}
	\sum_{n\ge 1}\frac{2\varphi((n,f))c_{N,f}(n)}{(n,f)n\varphi(4n)
	\prod_{\ell\mid (n,f)}\#C_N^{(\ell)}(1,1,f)}\\
&\quad\quad\ll\sum_{\substack{f\le mV\\ m\mid f\\ (f,2)=1}}\frac{\prod_{\ell\mid f}\#C_N^{(\ell)}(1,1,f)}{f^2\varphi(f)}
	\sum_{n>U}\frac{1}{\kappa_{2N}(n)\varphi(n)}\\
	&\quad\quad\quad+\sum_{\substack{f> mV\\ m\mid f\\ (f,2)=1}}
	\frac{\prod_{\ell\mid f}\#C_N^{(\ell)}(1,1,f)}{f^2\varphi(f)}
	\sum_{n\ge 1}\frac{1}{\kappa_{2N}(n)\varphi(n)},
\end{split}
\end{equation}
where the primes on the sums on $f$ are meant to indicate that the sums are to be restricted
to odd $f$ such that $\#C_N^{(\ell)}(1,1,f)\ne 0$ for all primes $\ell$ dividing $f$.

In~\cite[Lemma 3.4]{DP:1999}, we find that
\begin{equation*}
\sum_{n>U}\frac{1}{\kappa_1(n)\varphi(n)}
\sim\frac{c_0}{\sqrt U}
\end{equation*}
for some positive constant $c_0$.
In particular, this implies that the full sum converges.
From this we obtain a crude bound for the tail of the sum over $n$
\begin{equation*}
\begin{split}
\sum_{n>U}\frac{1}{\kappa_{2N}(n)\varphi(n)}
&=\sum_{\substack{kl>U\\ (l,2N)=1\\ \ell\mid k\Rightarrow\ell\mid 2N}}
	\frac{1}{\kappa_1(l)\varphi(l)\varphi(k)}
=\sum_{\substack{k\ge 1\\ \ell\mid k\Rightarrow\ell\mid 2N}}\frac{1}{\varphi(k)}
	\sum_{\substack{l>U/k\\ (l,2N)=1}}\frac{1}{\kappa_1(l)\varphi(l)}\\
&\ll\sum_{\substack{k\ge 1\\ \ell\mid k\Rightarrow\ell\mid 2N}}\frac{1}{\varphi(k)}
	\frac{\sqrt k}{\sqrt U}
\ll\frac{1}{\sqrt U}\prod_{\ell\mid N}\left(1+\frac{\ell}{(\ell-1)(\sqrt\ell-1)}\right)\\
&=\frac{1}{\sqrt U}\frac{N}{\varphi(N)}
	\prod_{\ell\mid N}\left(1+\frac{1}{\sqrt\ell(\ell-1)}\right)\left(1+\frac{1}{\sqrt\ell}\right)\\
&\ll\frac{1}{\sqrt U}\frac{N}{\varphi(N)}
	\prod_{\ell\mid N}\left(1+\frac{1}{\sqrt\ell}\right).
\end{split}
\end{equation*}
We have already noted that $N/\varphi(N)\ll\log\log N$.
It is a straightforward exercise as in~\cite[p.~63]{MV:2007} to show that
\begin{equation*}
\prod_{\ell\mid N}\left(1+\frac{1}{\sqrt\ell}\right)
<\exp\left\{O\left(\frac{\sqrt{\log N}}{\log\log N}\right)\right\}.
\end{equation*}
Thus, we conclude that
\begin{equation}\label{n tail}
\sum_{n>U}\frac{1}{\kappa_{2N}(n)\varphi(n)}
\ll\frac{N^\epsilon}{\sqrt U}
\end{equation}
for any $\epsilon>0$.
For the full sum over $n$, we need a sharper bound in the $N$-aspect, which we obtain
by writing
\begin{equation}\label{full n sum bound}
\begin{split}
\sum_{n\ge 1}\frac{1}{\kappa_{2N}(n)\varphi(n)}
&=\prod_{\ell\mid 2N}\left(1+\frac{\ell}{(\ell-1)^2}\right)
	\sum_{\substack{n\ge 1\\ (n,2N)=1}}\frac{1}{\kappa_{2N}(n)\varphi(n)}\\
&\le\frac{2N}{\varphi(2N)}\prod_{\ell\mid 2N}\left(1+\frac{1}{\ell(\ell-1)}\right)
	\sum_{n\ge 1}\frac{1}{\kappa_1(n)\varphi(n)}\\
&\ll\log\log N.	
\end{split}
\end{equation}

For any odd prime $\ell$ dividing $f$, we obtain the bounds
\begin{equation*}
\#C_N^{(\ell)}(1,1,f)
\le\begin{cases}
2\ell^{\nu_\ell(N)/2}&\text{if }\nu_\ell(f)>\nu_\ell(N)/2\text{ and }2\mid\nu_\ell(N),\\
\ell^{\nu_\ell(f)}&\text{otherwise}
\end{cases}
\end{equation*}
from equation~\eqref{valueofnumberC}.
Therefore, for every odd integer $f$, we have that
\begin{equation*}
\prod_{\ell\mid f}\#C_N^{(\ell)}(1,1,f)\le f,
\end{equation*}
and hence
\begin{equation}\label{f tail}
\sum_{\substack{f> mV\\ m\mid f\\ (f,2)=1}}\frac{\prod_{\ell\mid f}\#C_N^{(\ell)}(1,1,f)}{f^2\varphi(f)}
<\sum_{\substack{f>mV\\ m\mid f}}\frac{1}{f\varphi(f)}\ll\frac{1}{\varphi(m^2)V}.
\end{equation}
Substituting the bounds~\eqref{n tail},~\eqref{full n sum bound},
and~\eqref{f tail} into~\eqref{diff bn K and truncated K}, the lemma follows.
\end{proof}

\section{Concluding remarks.}\label{conclusion}

There are many open conjectures about the distribution of local invariants associated with the reductions modulo $p$
of a fixed elliptic curve defined over $\Q$  as $p$ varies over the primes.  Perhaps the most famous examples are the conjectures
of Koblitz~\cite{Kob:1988} and of Lang and Trotter~\cite{LT:1976}.
The Koblitz Conjecture concerns the number of primes $p\le X$ such that $\#E(\F_p)$ is prime.
The fixed trace Lang-Trotter Conjecture concerns the number of primes $p\le X$ such that
the trace of Frobenius $a_p(E)$ is equal to a fixed integer $t$.  Another
conjecture of Lang and Trotter (also called the Lang-Trotter Conjecture)
concerns the number of primes $p\le X$ such that the Frobenius field $\Q(\sqrt{a_p(E)^2-4p})$ is a fixed
imaginary quadratic field.

These conjectures are all completely open.  In order to gain evidence, it is natural to consider the averages for these
problems over some family of elliptic curves.
This has been done by various authors originating with the work of Fouvry and
Murty~\cite{FM:1996} for the number of supersingular primes (i.e., the fixed trace Lang-Trotter
Conjecture for $t=0$).  See~\cite{DP:1999,DP:2004,Jam:2004,BBIJ:2005,JS:2011,CFJKP:2011}
for other averages regarding the fixed trace Lang-Trotter Conjecture.
The average order for the Koblitz Conjecture was considered in~\cite{BCD:2011}.
Very recently, the average was successfully carried out for the Lang-Trotter Conjecture on Frobenius fields~\cite{CIJ:2013}.
The average order problems that we consider here and in~\cite{DS:2013} have very different features than
the above averages.
This is primarily because both $M_E(G)$ and $M_E(N)$ count finite sets of primes whose sizes vary with the parameters $G$ and $N$,
whereas the above problems seek estimates for the densities of (what are believed to be) infinite sets of primes.

Much like the conjectural constants appearing in the Twin Prime Conjecture and
the more general Bateman-Horn Conjectures~\cite{BH:1962}, the constants in the 
conjectural asymptotics for the
Lang-Trotter Conjectures and the Koblitz Conjecture can be written as Euler products in which the local Euler factors at each prime
$\ell$ can be understood in terms of the probability that $p$ 
satisfies the desired property modulo $\ell$.
Thus, for those conjectures, the constant from the average asymptotic gave strong evidence for the original conjectures as one
retrieves the local Euler factors of the conjectural constants (forgetting about a finite number of ``exceptional primes" for each elliptic curve).

In the questions considered in~\cite{DS:2013} and the present paper, it does not seem very likely that there should be an asymptotic for
a fixed elliptic curve.  Moreover, the ``constants" $K(N)$ and $K(G)$
(defined by~\eqref{defn of K(N)} and~\eqref{defn of K(G)}, respectively) are much more peculiar.
In the first place, they are not truly constant.
In the second place, it is not completely clear how to understand them in terms of local probabilities.
However, there are some parts of $K(N)$ which do seem to fit those local probabilities. 
More precisely,
given an elliptic curve $E/\Q$ without complex multiplication, for all but finitely many primes $\ell$, we have an isomorphism
$\mathrm{Gal}(\Q(E[\ell]))\isom\GL_2(\Z/\ell\Z),$
where $\Q(E[\ell])$ denotes the field obtained by adjoining to $\Q$ the coordinates of the $\ell$-division torsion points $E[\ell]$.
Furthermore, $\#E_p(\F_p)=p+1-a_p(E) \equiv \det(\sigma_p)+1-\mathrm{tr}(\sigma_p) \pmod \ell$,
where $\sigma_p$  denotes any element of the image of the Frobenius class at $p$ in $\GL_2(\Z/\ell\Z)$.
Letting
\begin{equation*}
C_N(\ell):=\{g\in\GL_2(\Z/\ell\Z): \det(g)+1-\mathrm{tr}(g)\equiv N\pmod\ell\},
\end{equation*}
one readily finds in the case that $\ell\nmid N$,
\begin{eqnarray*}
\frac{\prob ( p+1-a_p(E) \equiv N \pmod \ell )}{\prob( n \equiv N \mod \ell)} &=&
\frac{\displaystyle \frac{\#C_N(\ell)}{\#\GL_2(\Z/\ell\Z)}}{\displaystyle \frac{1}{\ell}}\\
&=& \left(1-\frac{\leg{N-1}{\ell}^2\ell+1}{(\ell-1)^2(\ell+1)}\right),
\end{eqnarray*}
where the denominator is the probability that a random integer $n$ is congruent to $N$ modulo $\ell$.
That is, we find agreement with the Euler factors of $K(N)$ at the primes $\ell$ not dividing $N$. 
It should be noted that Euler factors of $K(N)$ and $K(G)$ agree that the primes $\ell$ not dividing $N=\#G$.
The local factors of $K(N)$ at the primes $\ell$ dividing $N$ are more subtle. 
Some of the issues appearing in the analysis of those local factors are addressed in ~\cite{MPS},
where a statistical analysis of the function $K(N)N/\varphi(N)$ has been carried out.
Employing a new technique to address moments of functions which are ``almost but not quite multiplicative," the authors show that
\begin{equation*}
\sum_{\substack{N \leq x\\ 2\nmid N}}K(N)\frac{N}{\varphi(N)\log N}\sim\frac{1}{3}\frac{x}{\log x}.
\end{equation*}
Note that this is what one should expect as exactly $1/3$ of the elements of $\GL_2(\Z/2\Z)$ have odd trace
(i.e., satisfy the condition $\det(g)+1-\mathrm{tr}(g)\equiv 1\pmod 2$). In collaboration with
the authors of \cite{MPS}, the authors of the present paper are studying further properties of
the function $K(N)$ which could lead to a probabilistic model for the local factor at the primes
$\ell$ dividing $N$. It seems that those could eventually be understood by looking at probabilities
on the group $\GL_2(\Z_\ell)$.

As far as we know, no statistical analysis similar to~\cite{MPS} has been made for $K(G)$,
but this could certainly be an interesting avenue of research.

\bibliographystyle{alpha}
\bibliography{references}


\end{document}